\documentclass[final,intlimits]{article}
\usepackage{amsmath, amssymb}

\title{\bf POSITIVE FORMS ON BANACH SPACES}
\author{B\'alint Farkas, M\'at\'e Matolcsi}
\date{25th June 2001}
\newtheorem{theo}{Theorem}
\newtheorem{lem}{Lemma}
\newtheorem{defi}{Definition}

\newcommand\dom{\mathrm{Dom}\:}
\newcommand\ran{\mathrm{Ran}\:}

\newcommand{\nam}{\textsc}
\newcommand{\EE}{\mathrm{E}\:}
\newcommand{\aaa}{\mathcal{A}}
\newcommand{\dd}{\:\mathrm{d}}
\newcommand{\adj}{^\ast}
\newcommand{\aadj}{^{\ast\ast}}
\newcommand{\la}{[}
\newcommand{\ra}{]}
\newcommand{\ls}{(}
\newcommand{\rs}{)}
\newcommand{\nm}{\|}
\newcommand{\sq}{^\frac{1}{2}}
\newcommand{\fplus}{\stackrel{.}{+}}
\newcommand{\rr}{\mathbb{R}}
\newcommand{\cc}{\mathbb{C}}

\newcommand{\inv}{^{-1}}
\newcommand{\suphat}{^{\widehat{}}}
\newcommand{\lm}{\rightarrow}
\newcommand{\LL}{\mathcal{L}}
\newcommand{\der}{\partial}
\newcommand{\qed}{
\ifmmode %
\else \leavevmode\unskip\penalty9999 \hbox{}\nobreak\hfill
\fi
\quad\hbox{\rule{1ex}{1ex}}}

\newenvironment{proof}[1][Proof]{\setlength{\parskip}{0pt plus 0.4pt}
\normalfont\upshape
\trivlist
\item[\hskip\labelsep{\itshape#1.}\ignorespaces]}{\qed\endtrivlist\setlength{\parskip}{0pt plus 1pt}}

\newenvironment{rem}[1][\itshape{Remark}]{\setlength{\parskip}{0pt plus 0.4pt}
\normalfont\upshape
\trivlist
\item[\hskip\labelsep{#1.}\ignorespaces]}{\endtrivlist\setlength{\parskip}{0pt plus 1pt}}

\begin{document}
\maketitle

\begin{abstract}
The first representation theorem establishes a correspondence between 
positive, self-adjoint operators and closed, positive forms on Hilbert 
spaces. The aim of this paper is to show that some of the results 
remain true if the underlying space is a reflexive Banach space. In 
particular, the construction of the Friedrichs extension and the form 
sum of positive operators can be carried over to this case.
\end{abstract} 
\section{Introduction}
\label{sec:1}
Let $X$ denote a reflexive complex Banach space, and $X\adj$ its 
conjugate dual space (i.e. the space of all continuous, conjugate linear 
functionals over $X$). We will use the notation $\ls v,x):=v(x)$ for $v\in X\adj$ 
, $x\in X$, and $\ls x,v):=\overline{v(x)}$.  Let $A$ be a densely defined linear operator 
from $X$ to $X\adj$. Notice that in this context it makes sense to speak 
about positivity and self-adjointness of $A$. Indeed, $A$  defines a 
sesquilinear form on $\dom A \times \dom A$ via 
$t_{A}(x,y) = (Ax)(y) =\ls Ax, y)$ and $A$ is called positive if 
$t_{A}$ is positive, i.e. if $\ls Ax, x)\ge 0$ for all $x\in{\dom A}$. 
Also, the adjoint $A\adj$ of $A$ is defined (because $A$ is densely 
defined) and is a mapping from $X\aadj$ to $X\adj$, i.e. from $X$ to 
$X\adj$. Thus, $A$ is called self-adjoint if $A=A\adj$. Similarly, the 
operator $A$ is called symmetric if the form $t_{A}$ is symmetric.

In Section \ref{sec:2} we deal with closed, positive forms and associated 
operators, and we establish a generalized version of the first 
representation theorem. In Section \ref{sec:3} we apply the representation 
theorem in two situations: first we construct the Friedrichs extension 
of a positive, symmetric operartor, then we define the form sum of two 
positive, self-adjoint operators. With the help of a factorization lemma we give 
a definition of a partial order of positive, 
self-adjoint operators. The Friedrichs extension turns out to 
be the largest of all self-adjoint extensions with respect to this partial order.  
The factorization lemma also enables us to prove a certain commutation property of the 
form sum of positive self-adjoint operators.
In the last section we give two applications 
of the results, one in the theory of partial differential equations and one in probability theory.
Some results of this paper have already appeared in a not widely known paper \cite{birman}.
\section{Representation theorem}
\label{sec:2}
Let $D\subseteq{X}$ be a dense subspace, and let $t: D\times D \to \cc$ be 
a sesquilinear form on $D$ (where $t$ is linear in the first variable 
and conjugate linear in the second). Assume that $t$ is positive with 
positive lower bound, i.e. $t(x, x)\ge \gamma \|x\|^{2} , \  \gamma >0$. 
Assume also that $t$ is "closed" in the sense that $(D, t( \cdot  , \cdot  ))=:H$ is 
a Hilbert space (i.e. it is complete). In this case, the injection 
$i:H\to X$ is continuous so $H$ can be regarded as a subspace of $X$. 
For brevity we will use the notation $\la \cdot  , \cdot  \ra$ for $t( \cdot  , \cdot  )$. An operator 
$A$ from $X$ to $X\adj$ can be associated to the form $t$ in a natural 
way: let $x\in{D}$ and take the functional $\la x,y\ra , \  y\in{D}$; if this 
functional is continuous in the norm of $X$ then there is an element 
$z$ in $X\adj$ for which $\la x,y\ra=z(y)=:\ls z,y\rs$, in this case, let $Ax:=z$. 

\begin{theo}\label{thm:1}
With notations as above the operator $A: X\to X\adj$ is a positive, 
self-adjoint operator.
\begin{proof}
Let $v\in X\adj$ be an arbitrary element. Now, $\ls v,x\rs  \  x\in D$ is a 
continuous, conjugate linear functional on $H$. Indeed, 
$$|\ls v, x\rs |\le \| v \| \  \| x \| \le \frac{1}{\sqrt{\gamma}}\la x\ra  \ \| v\| = K\la x\ra,$$
where $\la x\ra$ denotes the norm of $H$, i.e. $\la x\ra = \la x,x\ra^{1/2}$. Thus, by the 
theorem of Riesz we have an element $f\in H$ such that $\ls v,x\rs =\la f,x\ra$. Define
an operator $B$ from $X\adj$ to $X$ by $Bv:=f$. Then $B$ is defined everywhere on 
$X\adj$, and $B$ is positive and bounded with $\|B\| \le \frac{1}{\gamma}$.
Indeed, $\ls z,Bz\rs =\la Bz,Bz\ra=\la Bz\ra^2 \ge 0$, and 
$$\|Bz\|^2 \le \frac{1}{\gamma}\la Bz\ra^2 =\frac{1}{\gamma}\ls Bz,z\rs \le \frac{1}{\gamma}\|Bz\| \ \|z\|.$$
Hence, $B$ is a bounded, positive, self-adjoint operator. Furthermore, $B$ is injective.
To see this, suppose that $Bz=0$. Then $0=\la Bz,g\ra=\ls z,g\rs $ for every $g\in H$, 
and $H$ is dense in $X$ therefore $z=0$. This means that the inverse $B^{-1}$ 
exists and is a linear mapping from $X$ to $X\adj$. We will show that $A=B^{-1}$.
Let $x\in \dom A$, then $\la x,y\ra=\ls t,y\rs $ for some $t\in X\adj$ and $Ax=t$. Also, $\ls t,y\rs =\la Bt,y\ra$ so 
$Bt=x$, and hence $A\subseteq B^{-1}$. Conversely, if $x\in \dom B^{-1}$ then $x=Bz$ 
for some $z\in X\adj$ and $\la x,y\ra=\la Bz,y\ra=\ls z,y\rs $ is continuous in $y$ therefore $x\in \dom A$ and 
$Ax=z=B^{-1}z$, which proves that $B^{-1}\subseteq A$.    
To complete the proof we have the following lemma, which is well known in Hilbert spaces.
\end{proof}
\end{theo}
\begin{lem}
If $B:X\adj\to X$ is a bounded, injective, self-adjoint operator then $A:=B^{-1}$ is also a 
self-adjoint operator from $X$ to $X\adj$.
\begin{proof}
First we show that $\ran B$ is dense in $X$. Indeed, if for some $v\in X\adj$
we have $\ls Bz,v\rs =0$ for every $z\in X\adj$, then $\ls Bz,v\rs =\ls z,Bv\rs =0$ so $Bv=0$ and $v=0$.
Hence $A$ is densely defined. Also, $A$ is symmetric, because if $x\in \dom A$ 
then $x=Bz$ for some $z\in X\adj$ and $\ls Ax,x\rs =\ls z,Bz\rs \in \rr$. Thus $A\subseteq A\adj$. 
To see the reverse inclusion, let $y\in \dom A\adj$ and let $x=Bz$ run through the 
elements of $\dom A$. Then $\ls Ax,y\rs =\ls z,y\rs $ and also $$\ls Ax,y\rs =\ls x,A\adj y\rs =\ls Bz,A\adj y\rs =\ls z, BA\adj y\rs $$
which means that $y=BA\adj y$, so $y\in \dom A$.
\end{proof}
\end{lem}         
\begin{rem} The previous arguments can be carried out whenever $(X,Y)$ is a dual pair of locally convex, topological linear spaces. In this case, we replace the condition of the lower bound by the natural assumption that the injection $i$ introduced above is continuous.
\end{rem}
\section{The Friedrichs extension and the form sum}
\label{sec:3}
In this section we apply the representation theorem in two situations. First we construct the 
Friedrichs extension of a densely defined positive operator.
\begin{theo}
Let $a:X\to X\adj$ be a positive, densely defined operator with positive lower bound,
$\ls ax,x\rs \ge \gamma \|x\|^2 , \ \gamma>0$ for every $x\in \dom a$. Then $a$ admits a positve 
self-adjoint extension with the same lower bound.
\begin{proof}
The form $t_a(x,y):=\ls ax,y\rs $ defines a pre-Hilbert space on $\dom a$. Denote the completion 
of this space by $H$, and the arising inner product by $\la  \cdot , \cdot \ra$. The injection 
$i:\dom a \to X$ extends by continuity to $H$ and the 
extension will be denoted by $I_a$. We prove that $I_a$ is injective. Notice first that
$\la t,y\ra=\ls at,I_ay\rs $ for all $t\in \dom a$, $y\in H$. Indeed, take a sequence 
$y_n\in \dom a$ , $y_n \to y$ in $H$ (which implies convergence in $X$ as well), then
$$\ls at,I_ay\rs =\lim \ls at,I_ay_n\rs =\lim \la t,y_n\ra=\la t,y\ra.$$ Now assume that $I_ay=0$. Then 
$$\la y\ra^2=\lim \la y_n,y\ra=\lim \ls ay_n,I_ay\rs =0$$ therefore $y=0$ which means that $I_a$ is injective.
Thus $H$ can be regarded as a subspace of $X$ and Theorem \ref{thm:1} can be applied. 
It is clear that the arising self-adjoint operator $A_F$ is an extension of $a$ and we also see 
from the proof of Theorem \ref{thm:1} that $\ls A_Fx,x\rs \ge \gamma \|x\|^2$ for all $x\in \dom A$. This operator 
will be called the Friedrichs extension of $a$.
\end{proof}
\end{theo}
                    %\subsection{Partial order of positve operators}

Next we examine a characterizing property of the Friedrichs extension. In the Hilbert space 
setting the Friedrichs extension distinguishes itself by being the largest, self-adjoint 
extension with respect to the usual partial ordering of positive operators. 
The definition of this ordering is that for positive operators $A$ and $B$ we have $A\ge B$ 
if $\dom A^{\frac{1}{2}}\subseteq \dom B^{\frac{1}{2}}$ and 
$\ls A^{\frac{1}{2}}x,A^{\frac{1}{2}}\rs \ge \ls B^{\frac{1}{2}}x, B^{\frac{1}{2}}\rs $ for all 
$x\in \dom A^{\frac{1}{2}}$. We will now examine how this definition can be carried over 
to our situation. The following factorization lemma is well known in Hilbert spaces (see e.g. \cite{sebprok1}, \cite{sebstoch}). For bounded positive self-adjoint operators from $X$ to $X\adj$ this lemma was also 
proved in \cite{vakh}.  
and it plays a key role in the characterization of covariance operators 
of Banach space valued random variables.
\begin{lem}\label{lem:2}
Let $A$ be a positive self-adjoint operator from $X$ to $X\adj$ (that is \/$\ls Ax,x\rs \ge 0$ for all $x\in \dom A$). 
Then there exists an auxiliary Hilbert space $H$ and an operator $J:H\to X\adj$ such that $A=JJ\adj$.              
\begin{proof}
Define an inner product on $\ran A$ by $\la Ax,Ay\ra:=\ls Ax,y\rs $. It is well defined because if 
$Ax_1=Ax_2$ and $Ay_1=Ay_2$ then $$\ls Ax_1,y_1\rs =\ls Ax_2,y_1\rs =\ls x_2,Ay_1\rs =\ls x_2,Ay_2\rs =\ls Ax_2,y_2\rs .$$
Furthermore it is positive definite, because if $\la Ax,Ax\ra=\ls Ax,x\rs =0$ then by the Cauchy 
inequality we have $$|\ls Ax,y\rs |^2\le \ls Ax,x\rs \ls Ay,y\rs =0$$ for all $y\in \dom A$ which implies that 
$Ax=0$. Thus $(\ran A , \la  \cdot , \cdot \ra)$ is a pre-Hilbert space. Denote the 
completion of this space by $H_A$. Define the operator $J: H_A\to X\adj$ by $\dom J= \ran A$ and 
$J(Ax):=Ax$ for all $Ax\in \ran A$. Then, by definition 
$\dom J\adj=\{ y\in X : |\ls Ax,y\rs |^2\le M_y\ls Ax,x\rs   \mbox{ for all }  x\in \dom A \}$, 
in particular 
$\dom A\subseteq \dom J\adj$ and $J\adj y=Ay$ for all $y\in \dom A$. Thus $JJ\adj$ is an extension of $A$
and $JJ\adj$ is symmetric. It is also clear that a self-adjoint operator is maximal symmetric 
just as in the context of Hilbert spaces. This means that $A=JJ\adj$.
\end{proof}
\end{lem}
\begin{rem}\rm Notice that in the context of Hilbert spaces $\dom A^{\frac{1}{2}}=\dom J\adj$ 
and $\ls A^{\frac{1}{2}}x,A^{\frac{1}{2}}x\rs =\la J\adj x,J\adj x\ra$.
\end{rem}
Now we are in position to give a definition of ordering without the use of square roots.
\begin{defi}
For positive, self-adjoint operators $A$ and $B$ we say that $A\ge B$ if and only if 
$\dom J_A\adj\subseteq \dom J_B\adj$ and $\la J_A\adj y,J_A\adj y\ra_A\ge \la J_B\adj y,J_B\adj y\ra_B$ for all
$y\in \dom J_A\adj$.
\end{defi}
In order to understand this definition better it would be desirable to give a characterization 
of $\dom J_A\adj$ and $\la J_A\adj y,J_A\adj y\ra_A$ in terms of $A$ only.
\begin{lem}
\label{lem:3}
With notations as above we have 
$$\dom J_A\adj= \left\{ y\in X : \sup_{x\in \dom A, \ls Ax,x\rs \le 1} |\ls Ax,y\rs |^2<\infty \right\}$$ and 
$$\la J_A\adj y,J_A\adj y\ra_A=\sup_{x\in \dom A, \ls Ax,x\rs \le 1} |\ls Ax,y\rs |^2$$           
\begin{proof}The characterization of $\dom J\adj$ is clear from 
$$\dom J\adj=\{ y\in X : |\ls Ax,y\rs |^2\le M_y\ls Ax,x\rs  \mbox{ for all }  x\in \dom A \}.$$ 
To see the other equality notice that $\ran A$ is dense in $H_A$, therefore we have
$$\la J_A\adj y,J_A\adj y\ra_A=\sup_{\ls Ax,x\rs \le 1} |\la J\adj y,Ax\ra|_A^2=\sup_{\ls Ax,x\rs \le 1} |\ls y, Ax\rs |^2$$
\end{proof}
\end{lem}
\begin{rem}\rm Notice that this new definition of partial ordering coincides with the usual 
one when $X$ is a Hilbert space. If we allow the value of the indicated supremum to be $+\infty$,
then we can say in short that $A\ge B$ if and only if 
$$\sup_{x\in \dom A, \ls Ax,x\rs \le 1} |\ls Ax,y\rs |^2 \ge \sup_{x\in \dom B, \ls Bx,x\rs \le 1} |\ls Bx,y\rs |^2$$
for all $y\in X$. We say that the left hand side of the equation is the form of $A$ on $X$ 
while the right hand side is that  of $B$.
\end{rem}
Also, notice that reflexivity and transitivity of the introduced relation is clear from 
Lemma \ref{lem:3}. To see that the relation is antisymmetric, assume that $A$ and $B$ have the same forms 
on $X$. Let $x\in \dom A$ and $y\in \dom B$. Then 
$\ls Ax,y\rs=\la J_A\adj x,J_A\adj y\ra=\la J_B\adj x,J_B\adj y\ra=\ls x,By\rs$ which means that $B\subseteq A\adj=A$
and hence $A=B$.    

Now, we establish the maximality of the Friedrichs extension with respect to this partial
order.
\begin{theo}
Assume that $a$ is a positive operator from $X$ to $X\adj$ with 
$\ls ax,x\rs \ge \gamma \|x\|^2 , \gamma >0$.
Then the Friedrichs extension is the largest, positive, self-adjoint extension of $a$  
with respect to the partial order introduced above.
\begin{proof}
The completion of $\ran A_F$ with inner product $\la A_Fx,A_Fy\ra_F:=\ls A_Fx,y\rs $ will be denoted by 
$H_F$.
Notice first that $\ran a$ is dense in $H_F$. Indeed, take $A_Fx\in H_F$ and let $x_n$ converge 
to $x=BA_Fx\in H$ in the norm of $H$ (recall that $H$ is the completion of 
$(\dom a , \la  \cdot , \cdot \ra)$). Then $ax_n=A_Fx_n$ is Cauchy in $H_F$ therefore $x_n\to y$ for some 
$y\in H_F$. Now, for all $z\in \dom A_F$ we have $\la A_Fz,A_Fx\ra=\ls A_Fz,x\rs $ and 
$$\la A_Fz,y\ra=\lim \la A_Fz,A_Fx_n\ra=\lim \ls A_Fz,x_n\rs =\ls A_Fz,x\rs $$ hence $A_Fx=y$. (In the last equality 
we used the fact that $x_n\to x$ in the norm of $X$ too, due to the positive lower bound of $a$.)

Next we prove that $\dom a$ is a core for $J_F\adj$. Take an arbitrary element 
$(x;J\adj x)\in X\times H_F$. By the argument above there exists a sequence $ax_n$ such that 
$ax_n\to J\adj x$ in $H_F$. Now, $x_n$ is Cauchy in $H$ so there is an element $z\in H$ such 
that $x_n\to z$ in $H$ (and consequently $x_n\to z$ in $X$). Also, $J_F\adj$ is a closed operator
therefore $z\in \dom J_F\adj$ and $J_F\adj z=J_F\adj x$. Thus, for all $u\in \dom A_F$ we have
$\ls A_Fu,x\rs =\la A_Fu,J_F\adj x\ra_F=\la A_Fu,J_F\adj z\ra=\ls A_Fu,z\rs $ and $A_F$ is surjective, so this implies 
$z=x$. Therefore $(x_n;ax_n)\to (x;J\adj x)$ in $X\times H_F$.

Now, take an arbitrary extension $A$ of $a$ and compare it to $A_F$. Define an operator $V$ from 
$H_F$ to $H_A$ by $\dom V=\ran a$ and $V(ax)=ax$. Then $V$ preserves norms and extends uniquely
 to an isometry on the whole of 
$H_F$. This extension will still be denoted by $V$. Notice that $J_A\adj x=ax$ for all 
$x\in \dom a$ so we have the inclusion $V(J_F\adj |_{\dom a})\subseteq J_A\adj$. Now, 
$\dom a$ is a core for $J_F\adj$ and $J_A\adj$ is a closed operator therefore 
$VJ_F\adj \subseteq J_A\adj$. This means that $\dom J_F\adj \subseteq \dom J_A\adj$ and 
$\la J_F\adj x\ra_F=\la J_A\adj x\ra_A$ so the form of $A_F$ is a restriction of that of $A$.      
\end{proof}
\end{theo}
Next we turn to the construction of the form sum of two positive self-adjoint operators (cf. \cite{us}).

Here we need a more general notion of closed forms. A positive form 
$t: D\times D \to \cc$ will be 
called closed if whenever $x_n\subseteq D$ and $x_n\to x$ in $X$ and $t(x_n-x_m, x_n-x_m)\to 0$ 
then $x\in D$ and $t(x_n-x, x_n-x)\to 0$ (notice that when $t$ has positive lower bound then 
this definition agrees with the previous one). In this more general context the representation 
theorem is no longer available (unless $X$ is a Hilbert space), but the form sum construction 
can be carried out if both forms 
are closed and at least one of them has positive lower bound.
     
Assume that $A$ is a positive self-adjoint operator with positive lower bound, and $B$ is an 
operator associated to a closed form $t_B$. Assume also that 
$H_{A,B}:=\dom J_A\adj \cap \dom t_B$ is 
dense in $X$. Then it is easy to see that $(H_{A,B}, t_A+t_B)$ is complete, thus the 
representation theorem can be applied. The arising positive self-adjoint operator will be called 
the form sum of $A$ and $B$, and will be denoted by $A\fplus B$. Next we will examine a certain 
commutation property of the form sum with the help of a factorization similar to that of Lemma \ref{lem:2}.

First we give a factorization of the form sum of positive, self-adjoint operators 
$A,B$ which have positive lower bound.
Let $J:H_A\oplus H_B\rightarrow X\adj$ be the densely defined operator, given by 
$J(Ax\oplus By)=Ax+By$ for $x\in \dom A,y\in\dom B$.
It is obvious that $J\adj$ exists and is densely defined. Indeed, it is not difficult to see that
$$\dom J\adj=\dom J_A\adj\cap \dom J_B\adj$$
Let us compute $J\adj$ on $\dom A\cap \dom B$. Take $z$ from the common domain of $A$ and $B$, 
and let $x,y$ run through $\dom A$ and $\dom B$ respectively. For brevity we will use the 
unified notation $\la \cdot , \cdot \ra$ for the inner products of the Hilbert spaces $H_A$, $H_B$ and 
$H_A\oplus H_B$. Now
$$
\ls J(Ax\oplus By),z)=\ls Ax+By,z\rs=\la Ax,Az\ra+\la By,Bz\ra=\la Ax\oplus By,Az\oplus Bz\ra,
$$
consequently $J\adj z=Az\oplus Bz$. 
Now, if $x\in\dom A\cap\dom B$, then
$$
J\aadj J\adj x=J\aadj (Ax\oplus Bx)=J(Ax\oplus Bx)=Ax+Bx,
$$
which means that $J\aadj J\adj$ is an extension of $A+B$.

We claim that $J\aadj J\adj$ is nothing else than the form sum $A\fplus B$. To see this compute 
the form determined by the symmetric operator $J\aadj J\adj$. Denote by $J_A$ and $J_B$ the 
factorizing operators for $A,B$ appearing in Lemma \ref{lem:2}. Since $J\aadj J\adj$ is 
symmetric, it suffices to show that $ A\fplus B\subseteq J\aadj J\adj$. 
Take $y\in\dom A\fplus B$, then it is straightforward that $y\in\dom J\adj$. We should 
determine $\la J\adj y,J\adj y\ra$, the following calculation gives an upper estimate.
\begin{equation*}
\begin{split}
\la J\adj y&,J\adj y\ra=\\
\sup&\left\{|\la Au\oplus Bv,J\adj y\ra |^2:u\in\dom A, v\in\dom B,\ls Au,u\rs+\ls Bv,v\rs\leq 1\right\}=\\
\sup&\left\{|\ls Au+Bv, y\rs |^2:u\in\dom A, v\in\dom B,\ls Au,u\rs+\ls Bv,v\rs\leq 1\right\}\leq \\
\sup&\left\{|\la (J_A\adj\oplus J_B\adj)(u\oplus v), (J_A\adj\oplus J_B\adj)(y\oplus y)\ra |^2: u\in \dom J_A\adj , v\in\dom J_B\adj,\right.\\&\quad  \left.\la (J_A\adj\oplus J_B\adj)(u\oplus v)\ra^2\leq 1\right\}=\\
\la (J_A\adj&\oplus J_B\adj)(y\oplus y)\ra^2=\la J_A\adj y,J_A\adj y\ra+\la J_B\adj y,J_B\adj y\ra=\ls A\fplus B y,y\rs
\end{split}
\end{equation*}
For the reverse inequality set 
$\lambda=\la J_A\adj y\ra^2/(\la J_A\adj y\ra^2+\la J_B\adj y\ra^2)$ and write:
\begin{equation*}
\begin{split}
\la J\adj y&,J\adj y\ra=\\
\sup&\left\{|\la Au\oplus Bv,J\adj y\ra |^2:u\in\dom A, v\in\dom B,\ls Au,u\rs+\ls Bv,v\rs\leq 1\right\}\geq\\
\sup&\left\{|\la J_A\adj u,J_A\adj y\ra+\la J_B\adj v,J_B\adj y\ra|^2:\ls Au,u\rs\leq\lambda, \ls Bv,v\rs\leq 1-\lambda\right\}=\\
\sup&\left\{(|\la J_A\adj u,J_A\adj y\ra|+|\la J_B\adj v,J_B\adj y\ra|)^2:\ls Au,u\rs\leq\lambda, \ls Bv,v\rs\leq 1-\lambda\right\}=\\
\lambda \la J_A\adj& y\ra^2+(1-\lambda) \la J_B\adj y\ra^2+2\sqrt{\lambda(1-\lambda)}\la J_A\adj y\ra\la J_B\adj y\ra=\\\la J_A\adj& y\ra^2+\la J_B\adj y\ra^2=\ls A\fplus B y,y\rs
\end{split}
\end{equation*}
Combining these two estimates, we see that
$$
\la J\adj y,J\adj y\ra=\ls A\fplus B y,y\rs
$$ 
for all $y\in\dom A\fplus B$. This implies that $\la J\adj x,J\adj y\ra=\ls A\fplus B x,y\rs$ 
for all $x,y\in\dom A\fplus B$, and as consequence we have that 
$A\fplus B\subseteq J\aadj J\adj$. Thus we have obtained the following theorem.

\begin{theo} $J\aadj J\adj=A\fplus B$ and it is an extension of $A+B$.
\end{theo}

Now, we observe that commutation with bounded operators is preserved when taking the form sum 
of two positive, self-adjoint operators. Before all, we characterize this commutation property 
by means of a bounded operator on an auxiliary Hilbert space.

Let $A:X\to X\adj $ be a positive, selfadjoint operator 
(in the next two lemmas we do not require $A$ to have positive lower bound). 
%                        with lower bound $\gamma>0$. 
Suppose that $E$ is a bounded, linear operator on $X$ which leaves $\dom A$ invariant and 
satisfies
\begin{equation}
\label{eq:7}
E\adj A\subseteq A E.
\end{equation}
In this case, define $\hat{E}$ on $\ran A \subset H_A$ as $\hat{E}(Ax):=AEx$. This is indeed a 
definition of a  
linear operator, as will be seen in the following lemma.
\begin{lem} $\hat{E}$ is well defined and furthermore is a bounded linear operator on $\ran A$.
\begin{proof}
Fix $x\in\dom A$ and consider the following:
\begin{equation*}
\label{eq:5}
\begin{gathered}
\la \hat{E}Ax, \hat{E}Ax\ra=
\la AEx,AEx\ra=
\ls AEx,Ex\rs=\ls E\adj Ax,Ex\rs=\ls Ax,E^2x\rs=\\
\la Ax,AE^2x\ra\leq
\la Ax,Ax\ra\sq\la AE^2x,AE^2x\ra\sq=
%                            \la Ax,Ax,\ra\sq\la \hat{E}^2 Ax,\hat{E}^2Ax\ra\sq
\end{gathered}
\end{equation*}
By induction, we conclude the following:
\begin{equation*}
\label{eq:6}
\begin{gathered}
\la \hat{E} Ax,\hat{E} Ax\ra\leq \la Ax,Ax\ra^{1-\frac{1}{2^n}}\la AE^{2^n}x, AE^{2^n}x\ra^{\frac{1}{2^n}}\leq\\\la Ax,Ax\ra^{1-\frac{1}{2^n}}\ls Ax,E^{2^{n+1}}x\rs^{\frac{1}{2^n}}\leq\la Ax,Ax\ra^{1-\frac{1}{2^n}}\nm Ax\nm^{\frac{1}{2^n}}\nm E^{2^{n+1}} x\nm^{\frac{1}{2^n}}
\end{gathered}
\end{equation*}
Now, taking limit in $n$ we obtain 
$$\la \hat{E} Ax,\hat{E} Ax\ra\leq \la Ax\ra^2 r(E^2),$$
where $r(E^2)$ denotes the spectral radius of $E$. This implies that $\hat{E}$ is well-defined 
and continuous on $\ran A$.
\end{proof}
\end{lem}

Extend $\hat{E}$ to $H_A$, and denote this bounded operator also by $\hat{E}$. 
We connect the self-adjointness of $\hat{E}$ with the commutation property \eqref{eq:7} of $E$ 
in the following (cf. \cite{sebestyen}).

\begin{lem} $\hat{E}$ is a self-adjoint operator in $H_A$.
\begin{proof}
It suffices to check that $\hat{E}Ax=\hat{E}\adj Ax$. So let $x\in\dom A$ be fixed and 
$y\in\dom A$ arbitrary, then
$$
\la Ax,\hat{E}Ay\ra=\la Ax,AEy\ra=\ls Ax,Ey\rs=\ls E\adj Ax,y\rs=\ls AEx,y\rs=\la \hat{E}Ax, Ay\ra,
$$ 
implying the desired equality.
\end{proof}
\end{lem}

With the help of the previous arguments, we obtain the theorem on the commutation property of 
the form sum.

\begin{theo} Let $A,B$ be two positive, self-adjoint operators with positive lower bound, and 
assume that a bounded, linear operator $E$ is given on $X$, such that \eqref{eq:7} is satisfied 
together with a similar condition concerning $B$. Then for the form sum $A\fplus B$, we have 
$$E\adj (A\fplus B)\subseteq(A\fplus B) E.$$
\begin{proof}
Denote the bounded, self-adjoint operators resulting from the above construction on $H_A$ and 
$H_B$ by $\hat{E}_A$ and $\hat{E}_B$ respectively. Consider the bounded, self-adjoint operator 
$\hat{E}_A\oplus\hat{E}_B$ on the Hilbert space $H_A\oplus H_B$. 
Then $E\adj J\subseteq J(\hat{E}_A\oplus \hat{E}_B)$ and 
$(\hat{E}_A\oplus \hat{E}_B)J\adj\subseteq J\adj E$, since for any $x\in \dom A,y\in \dom B$
\begin{equation*}
\begin{gathered}
J\left(\hat{E}_A\oplus \hat{E}_B\right)(Ax\oplus By)=J(AEx\oplus BEy)=AEx+BEy=\\E\adj Ax+E\adj By=E\adj J(Ax\oplus By),
\end{gathered}
\end{equation*}
and from this, using the continuity of $E\adj$, it follows that
$$
\left(\hat{E}_A\oplus \hat{E}_B\right)J\adj\subseteq\left[J\left(\hat{E}_A\oplus \hat{E}_B\right)\right]\adj
\subseteq  \left(E\adj J\right)\adj=J\adj E.
$$
We use these two inclusions and the boundedness of the self-adjoint operator 
$\hat{E}_A\oplus \hat{E}_B$ to conclude
\begin{equation*}
\begin{gathered}
E\adj J\aadj\subseteq\left(J\adj E\right)\adj\subseteq  \left[\left(\hat{E}_A\oplus \hat{E}_B\right)J\adj\right]\adj=J\aadj\left(\hat{E}_A\oplus \hat{E}_B\right)\mbox{ and }\\
E\adj\left(A\fplus B\right)=E\adj J\aadj J\adj\subseteq J\aadj \left(\hat{E}_A\oplus \hat{E}_B\right) J\adj\subseteq J\aadj J\adj E=\left(A\fplus B\right)E.
\end{gathered}
\end{equation*}
\end{proof}
\end{theo}  
Finally we observe the relation between the spectra of $\hat{E}=\hat{E}_A$ and that of $E$. 
(Here again the requirement that $A$ has psoitive lower bound is not necessary.) 
\begin{theo} The spectrum $\sigma(\hat{E})$ is contained in $\sigma(E)\cap\rr$. 
%                        while for the 
%                        point spectrum $P\sigma(\hat{E})$ we have the inclusion 
% $P\sigma(E\rest_{\dom A})\subseteq P\sigma(\hat{E})$.
\begin{proof}
Since $\hat{E}$ is self-adjoint it is clear that $\sigma(\hat{E})\subseteq\rr$. On the other 
hand, take any real $\lambda$  from the resolvent set of $E$ and $x\in\dom A$, then
$$
A(E-\lambda I)x=(E-\lambda I)\adj Ax,\mbox{ hence } A(E-\lambda I)\inv x=\left[\left(E-\lambda I\right)\inv\right]\adj Ax
$$
which means that we can define the operator $[(E-\lambda I)\inv]\suphat$. A short computation 
gives that
$$
\left[\left(E-\lambda I\right)\inv\right]\suphat=\left(\hat{E}-\lambda \hat{I}\right)\inv,
$$
indeed for $x\in\dom A$
\begin{equation*}
\begin{gathered}
\left[\left(E-\lambda I\right)\inv\right]\suphat \left(\hat{E}-\lambda \hat{I}\right)Ax=A\left(E-\lambda I\right)\inv\left(E-\lambda I\right)x=Ax,\mbox{ and}\\
\left(\hat{E}-\lambda \hat{I}\right)\left[\left(E-\lambda I\right)\inv\right]\suphat Ax=A\left(E-\lambda I\right)\left(E-\lambda I\right)\inv x=Ax.
\end{gathered}
\end{equation*}
This proves the statement.
%                So we have proved the first statement, 
%                while the second follows immediately from the definition 
%                 of $\hat{E}$.
\end{proof}
\end{theo}

\section{Application of the results}
\label{sec:4}
\subsection*{Covariance operators.} Consider a probability measure space 
$\langle \Omega,\aaa,\mu\rangle$, and let
$\xi:\Omega\rightarrow X$ a random variable i.e. a weakly measurable function. 
Suppose that $\xi$ possesses a weak expectation,
 in other words
 $$\EE\xi:=\int_\Omega \xi\dd\mu$$ exists as a Pettis integral. Note that if $X$ is reflexive, 
according to Dunford and Gelfand, this is equivalent
 to requiring the existence of
 $$\int_\Omega f(\xi)\dd\mu$$ for all $f\in X\adj$. Further, we make assumptions on the second 
moments,
 so suppose that the set
 $$D=\left\{f:f\in X\adj,\int_\Omega |f(\xi)|^2\dd\mu<+\infty\right\}
 $$
 is dense in $X\adj$.

 As an example, take $X=\ell_2$, $\Omega= \left\{\omega_n:n=1,2,\dots\right\}$ and 
$\mu(\{\omega_n\})=ce^{-(3/2)n}$ with a suitable constant $c$.
 Setting $\xi (\omega_n)_k=n^k/k!$, it is easy to compute that, in this case, $D\neq X\adj$ is 
dense.

 In the sequel we assume that $\EE\xi =0$, since we could take $\xi-\EE\xi$ instead of $\xi$.
 Define the sesquilinear form
 $$t(f,g)=\EE \left(f(\xi)\bar{g}(\xi) \right)$$
 for $f,g\in D$.
 \begin{theo} $t$ is a positive, closed, sesquilinear form on $D\times D$.
 \begin{proof}
 Positivity is trivial.
Suppose that $f_n\in D$ converges to $f\in X\adj$ and 
$\EE \left|f_n(\xi)-f_m(\xi)\right|^2\lm 0$, then $f_n(\xi)$ has a limit
$g\in \LL_2(\Omega,\mu)$, and  moreover $g$ and $f(\xi)$ conincide almost everywhere, 
hence $\EE |f(\xi)|^2<+\infty$, implying
$f\in D$ and $\EE \left|f_n(\xi)-f(\xi)\right|^2\lm 0$.
 \end{proof}
 \end{theo}
If $t$ possesses a positive lower bound and $X$ is reflexive, then the application of 
Theorem \ref{thm:1} provides a 
representing, self-adjoint operator $A$ from $X\adj$ to $X\aadj =X$, which is
called the covariance operator of the random variable $\xi$ (cf. \cite{vakh}). Note that if $X$ is a Hilbert 
space then the original version of the representation theorem provides the covariance operator 
of $\xi$ associated to the closed form $t$ (even if $t$ has lower bound 0).
%      requiring the positive
%      lower bound is not really a restricting condition, since we can add  
%      the original inner product 
%      to $t$ resulting a sesquilinear form, which is
%      strictly bounded from below.

One expects a characterization of the covariance operator of the sum of independent random 
variables. In view of our
result the following theorem is quite straightforward.
\begin{theo}
If $\xi$ and $\eta$ are independent random variables with covariance operators $A$ and $B$ 
respectively. Then the covariance operator of $\xi+\eta$ is
$A\fplus B$.
\end{theo}

\subsection*{Elliptic operators.}
This is a classical application of the Friedrichs extension (see e.g.
\cite{birman}). Take 
$X=L_p(\Omega ) , 1\leq p<+\infty$
where $\Omega$ is a bounded domain with smooth boundary in $\rr^n$. Define
the operator $A$ from 
$L_p(\Omega )$ to $L_q(\Omega )$ by $\dom A=C_0^\infty (\Omega)$ and 
$$Af=- \sum_{i,k=1}^n \frac{\der}{\der x_i} (a_{ik} \frac{\der
f}{\der x_k}) +bf$$ where 
$a_{ik}\in C^1(\Omega)$, $b\in L^1_{loc}(\Omega)$, $b\geq 0$ and
$$\sum_{i,k=1}^n a_{ik}(x)\beta_i \overline{\beta_k} \ge \gamma \sum_i^n
|\beta_i|^2 , \gamma >0$$
everywhere in $\Omega$ (uniform ellipticity). 
In this case we have
\begin{equation*}
\begin{split}
\ls Af,f\rs &=\int_\Omega \left(- \sum_{i,k=1}^n \frac{\der}{\der x_i} \left(a_{ik}
\frac{\der f}{\der x_k}\right) +bf\right)\overline{f} \dd x \\&=
\int_\Omega \left(\sum_{i,k=1}^n a_{ik} \frac{\der f}{\der
x_i}\frac{\overline{\der f}}{\der x_k} +b|f|^2 \right)\dd x \ge
\gamma \int_\Omega \sum_{i=1}^n \left|\frac{\der f}{\der x_i}\right|^2 \dd x.
\end{split}
\end{equation*}
Now, for $p\le 2n/(n-2)$ we have 
$$\int_\Omega \sum_{i=1}^n \left|\frac{\der f}{\der x_i}\right|^2 \dd x\ge c\|
f\|_p^2, \quad c>0 $$ 
by the Sobolev imbedding theorem (see e.g. \cite{adams} pp. 95-99).
Thus $A$ has positive lower bound. The Friedrichs extension of $A$ is
surjective, and this 
means that the equation 
$$- \sum_{i,k=1}^n \frac{\der}{\der x_i} \left(a_{ik} \frac{\der f}{\der
x_k}\right) +bf=g$$ 
has a weak solution for every $g\in L_q(\Omega)$ whenever $q\ge 2n/(n+2)$.

%                        \subsection*{Scr\"{o}dinger operators.}

\vspace{1cm}
\noindent
\nam{B\'alint Farkas}\\
\small\noindent\nam{Department of Applied Analysis, E\"otv\"os Lor\'and University}\\
\small\noindent {Kecskem\'eti u. 10-12, 1053 Budapest, Hungary}\\
\small\noindent e-mail: \texttt{fbalint@cs.elte.hu}\\
\vspace{1em}

\normalfont
\noindent\nam{M\'at\'e Matolcsi}\\
\small\noindent\nam{Department of Applied Analysis, E\"otv\"os Lor\'and University}\\
\small\noindent {Kecskem\'eti u. 10-12, 1053 Budapest, Hungary}\\
\small\noindent e-mail: \texttt{matomate@cs.elte.hu}
\end{document}